\documentclass{amsart}
\usepackage{amsmath}
\usepackage{amsthm}
\usepackage{amssymb}
\usepackage{graphics}
\usepackage{psfig}

\DeclareMathOperator{\supp}{supp}
\theoremstyle{plain}
  \newtheorem{theorem}{\bf Theorem}[section]
  \newtheorem{proposition}[theorem]{\bf Proposition}
  \newtheorem{lemma}[theorem]{\bf Lemma}

\theoremstyle{remark}

\begin{document}
\title{Some special solutions of the Schr\"odinger equation}
\author{J. A. Barcel\'o, J. M. Bennett, A. Carbery, A. Ruiz and M. C. Vilela}\thanks{
All authors were supported by the EC project ``HARP''. The first
and fourth were also supported by Spanish Grant BFM2002-02204, the
fifth by Spanish Grant MTM2004-03029, the second by EPSRC
Postdoctoral Fellowship GR/S27009/02 and the third by a Leverhulme
Study Abroad Fellowship and EC project ``Pythagoras II".}
\address{J.A. Barcel\'o, ETSI de Caminos, Universidad Polit\'ecnica de
Madrid, 28040, Madrid, Spain.} \email{juanantonio.barcelo@upm.es}
\address{J. M. Bennett, School of Mathematics, The Watson Building, University of Birmingham, Edgbaston,
Birmingham, B15 2TT, England.} \email{J.Bennett@bham.ac.uk}
\address{A. Carbery, School of Mathematics, University of Edinburgh,
JCMB, King's Buildings, Mayfield Road, Edinburgh, EH9 3JZ,
Scotland.} \email{A.Carbery@ed.ac.uk}
\address{A. Ruiz, Departamento de Matem\'aticas, Universidad Aut\'onoma
de Madrid, 28049, Madrid, Spain.}\email{alberto.ruiz@uam.es}
\address{M. C. Vilela, Departamento de Matem\'atica Aplicada, Universidad de Valladolid, Plaza Santa Eulalia 9 y 11, Segovia, Spain.}
\email{maricruz@dali.eis.uva.es}
\date{1st April 2006}
\begin{abstract}
We describe a certain ``self-similar" family of solutions to the
free Schr\"odinger equation in all dimensions, and derive some
consequences of such solutions for two specific problems.
\end{abstract}
\maketitle
\section{Introduction}
In this paper we describe certain ``self-similar" solutions to the
initial value  problem associated to the free Schr\"odinger
equation:
\begin{equation}\label{ecuacion}
\begin{cases}
i \partial_t u+\Delta_x u=0 \hspace{0.7cm} (x,t)\in \mathbb
{R}^{n-1}
\times \mathbb {R}\hspace{0.7cm} n \geq 2, \\
u(x,0)=f(x).
\end{cases}
\end{equation}
A straightforward application of the Fourier transform allows us
express the solution $u$ of \eqref{ecuacion} as
\begin{equation*}
u(x,t)= \int_{{\mathbb{R}}^{n-1}} e^{- \pi i t |\xi|^2 + 2 \pi i
x\cdot\xi} \hat{f}(\xi) d \xi,
\end{equation*}
where $\hat{f}$ is the Fourier transform of $f$. As usual, we
denote this solution by $e^{it \Delta}f(x)$.

In Section 2 we describe in detail the particular solutions that
we have in mind, and in Section 3 we give applications to two
established problems in harmonic analysis and the theory of
Schr\"odinger equations. These problems concern weighted estimates
for solutions to \eqref{ecuacion} in two variants, and are related
to the restriction/extension operators for the base of the
paraboloid.


\subsubsection*{Notation} For non-negative quantities $X$ and $Y$
we use $X\lesssim Y$ ($X\gtrsim Y$) to denote the existence of a
positive constant $C$, depending on at most $n$, such that $X\leq
CY$ ($X\geq CY$). We write $X\sim Y$ if both $X\lesssim Y$ and
$X\gtrsim Y$.

\section{Special solutions}
\label{especiales}
Let $0<\delta<<1$ and $0<\sigma<\frac{1}{2}$. We consider the
function of one variable
\begin{equation}\label{alberto}
g= \sum_{\ell \in \mathbb{N},\;1 \leq \ell \leq \delta^{- \sigma}
}
\chi_{_{(\ell\delta^{\sigma}-\delta,\ell\delta^{\sigma}+\delta)}}.
\end{equation}
Note that $g$ is simply the characteristic function of a union of
disjoint, equally spaced subintervals of $[0,1]$ of equal size. We
take as initial data in \eqref{ecuacion} $f$ such that
\begin{equation}\label{exs}
\hat{f}(\xi)= \prod_{j=1}^{n-1} g(\xi_j),
\end{equation}
where $\xi=(\xi_1,\hdots,\xi_{n-1})$. The corresponding solution
of the free Schr\"odinger equation is given by
\begin{equation}\label{1.2}
e^{it \Delta}f(x)= \prod_{j=1}^{n-1} \int_{0}^{1} e^{- \pi i t
\xi_j^2 + 2 \pi i x_j.\xi_j} g(\xi_j) d \xi_j.
\end{equation}
We wish to identify a set
$\Omega\subset\mathbb{R}^{n-1}\times\mathbb{R}$ upon which
$|e^{it \Delta}f(x)|$ is ``large". In order to achieve this we
look for points $(x,t)$ for which there is essentially no
cancellation in the above integrals; i.e. for which the phases
$-t\xi_j^2/2+x_j\xi_j$ are within a small (say $1/10$)
neighbourhood of $\mathbb{Z}$ for all $\xi_j\in\supp(g)$ and
$1\leq j\leq n-1$. By the product structure of \eqref{1.2}, it
suffices to consider the integral
$$
\int_{0}^{1} e^{- \pi i t \xi^2 + 2 \pi i s\xi} g(\xi) d\xi,
$$
where $s\in\mathbb{R}$.

If $\xi\in\supp(g)$, then $\xi = \ell\delta^{\sigma} +
\varepsilon$ for some positive integer $\ell\leq\delta^{- \sigma}$
and $ |\varepsilon| \leq \delta$. Now consider $s,t\in\mathbb{R}$
of the form $s=p\delta^{-\sigma}$ and $ t=2q\delta^{-2\sigma}$,
where $p,q\in\mathbb{N}$. For such $s$, $t$  and $\xi$ we have
$$
s\xi-t\xi^2/2= p\ell+p\delta^{-\sigma}\varepsilon-
(q\ell^2+2q\ell\delta^{-\sigma}\varepsilon+q\delta^{-2\sigma}\varepsilon^2).$$
Since $p\ell, q\ell^2\in\mathbb{N}$, we are interested in the
values of $p$ and $q$ for which
\begin{equation}\label{juan}
|p\delta^{-\sigma}\varepsilon|,
|q\ell\delta^{-\sigma}\varepsilon|,
|q\delta^{-2\sigma}\varepsilon^2|\leq c\;\mbox{ for all }\; 1\leq\
\ell\leq\delta^{-\sigma}\;\mbox{ and } \;|\varepsilon|<\delta,
\end{equation}
where $c$ is a positive constant, such as $1/40$ say. As is easily
verified, \eqref{juan} holds precisely when
$|p|\lesssim\delta^{\sigma-1}$ and
$|q|\lesssim\delta^{2\sigma-1}$, and hence for such $s$ and $t$,
\begin{equation}\label{jon}
\left|\int_{0}^{1} e^{- \pi i t \xi^2 + 2 \pi i s\xi} g(\xi)
d\xi\right|\sim\delta^{1-\sigma}.
\end{equation}

\begin{figure}[h]
{\centerline{\psfig{figure=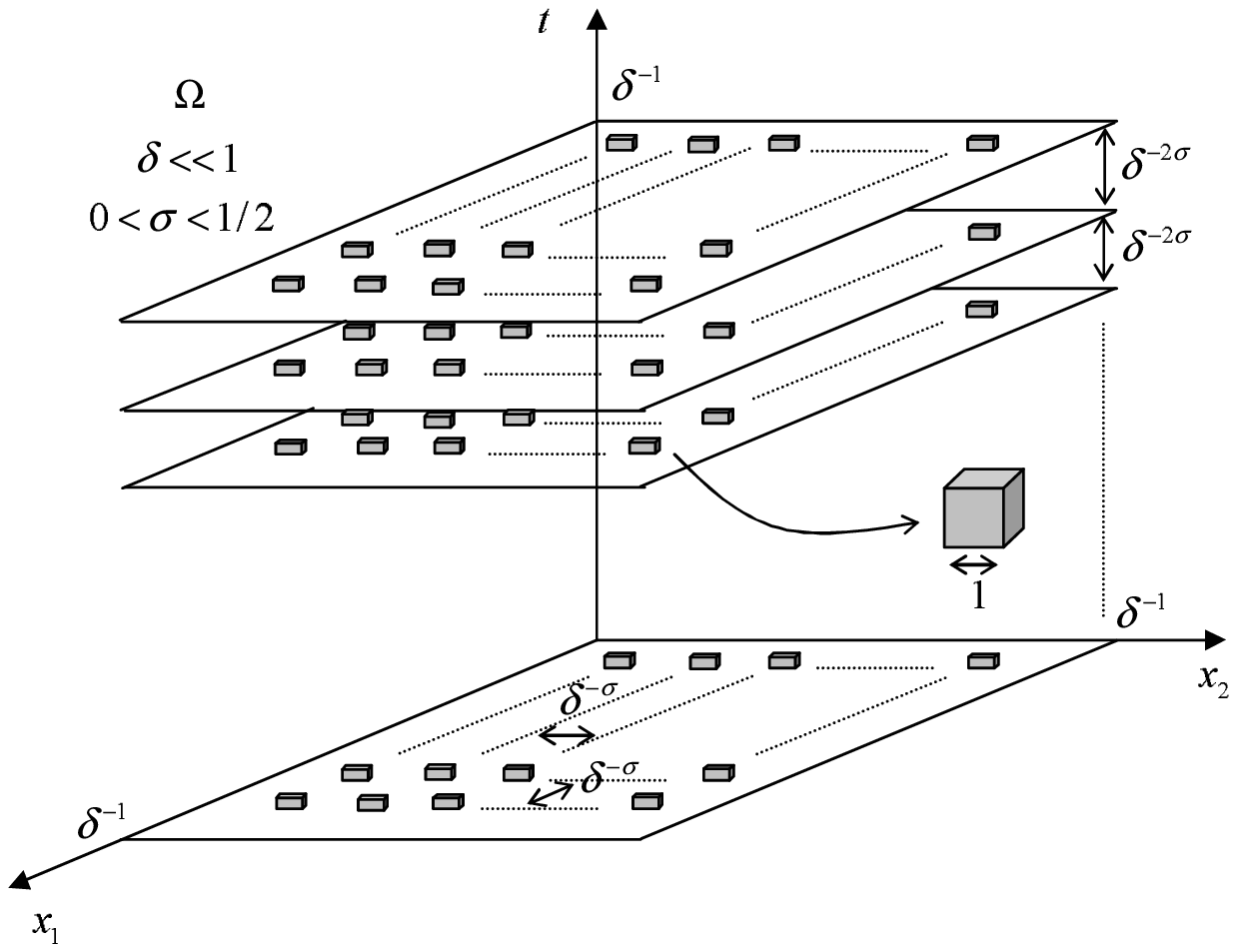,height=9.5cm,width=12.5cm}}}
\caption{The set $\Omega$ is the union of
$\delta^{2\sigma-1}\delta^{(\sigma-1)(n-1)}$ balls of radius $1$
painted in grey.} \label{Omega}
\end{figure}

We now define
$X=\{p\delta^{-\sigma}:p\in\mathbb{N}\mbox{ with } p\lesssim\delta^{\sigma-1}\}$, and
$$\Lambda= \{ (x,t)\in
\mathbb{R}^{n-1}\times\mathbb{R}: x\in X^{n-1}\mbox{ and } t = 2q
\delta^{-2 \sigma} \mbox{ where }q \in \mathbb{N} \mbox{ and }
q\lesssim\delta^{2\sigma -1}\}.
$$
Hence by \eqref{jon}, $ |e^{it
\Delta}f(x)|\sim\delta^{(n-1)(1-\sigma)}, $ for all
$(x,t)\in\Lambda$. Now since $e^{it \Delta}f(x)$ may be viewed as
the Fourier transform of a certain compactly supported
measure\footnote{See \eqref{red} for an explicit expression of
this.} on $\mathbb{R}^{n}$, this estimate continues to hold for
$(x,t)$ belonging to an $O(1)$ neighbourhood of $\Lambda$. We
denote this union of $O(1)$-balls by $\Omega$ (see Figure
\ref{Omega}).

We notice that for $t$ belonging to the projection of $\Omega$ onto the
$t$-axis, $\|e ^{it \Delta} f\|_2 \sim \|f\|_2$, so that the mass of
$ |e^{it \Delta}f(\cdot)|^2 $ on the section through $\Omega$ at height $t$
accounts for a positive proportion of its total mass at height $t$.

We now consider a succession of such examples of increasing
complexity.

Suppose now that $\delta>0$ is such that $1/\delta$ is an integer,
and let $k\in\mathbb{N}$ be fixed. (We allow the implicit
constants in the $\lesssim$ $\gtrsim$ $\sim$ notation to depend on
$k$.) Consider the function of one variable
\begin{equation}\label{tony}
g_k= \sum_{1\leq \ell_1,\hdots,\ell_{k} \leq \delta^{- \sigma} }
\chi_{\{s:|s-
(\ell_1\delta^{\sigma}+\ell_2\delta^{\sigma+1}+\cdots+\ell_k\delta^{\sigma+k-1})|
\leq \delta^{k}\}}.
\end{equation}
Note that $g_k$ is ``self-similar" in the sense that for $k\geq
2$, $g_k$ may be written as a certain linear combination of
rescaled and translated copies of $g_{k-1}$, where $g_1$ coincides
with the function $g$ defined in \eqref{alberto}.

As before we take as initial data in \eqref{ecuacion} $f_k$ such
that
\begin{equation*}
\hat{f_k}(\xi)= \prod_{j=1}^{n-1} g_k(\xi_j),
\end{equation*}
where $\xi=(\xi_1,\hdots,\xi_{n-1}).$ The corresponding solution
of the free Schr\"odinger equation is given by
\begin{equation}\label{1.3}
e^{it \Delta}f_k(x)= \prod_{j=1}^{n-1} \int_{0}^{1} e^{- \pi i t
\xi_j^2 + 2 \pi i x_j.\xi_j} g_k(\xi_j) d \xi_j.
\end{equation}

As before it suffices to consider the integral
$$
\int_{0}^{1} e^{- \pi i t \xi^2 + 2 \pi i s\xi} g_k(\xi) d\xi,
$$
where $s\in\mathbb{R}$.

If $\xi\in\supp(g)$, then
$$
\xi = \sum_{r=1}^k\ell_r\delta^{\sigma+r-1} + \varepsilon$$ for
some positive integers $\ell_r\leq\delta^{- \sigma}$ and $
|\varepsilon| \leq \delta^k$. Now consider $s,t\in\mathbb{R}$ of
the form
$$
s=\sum_{m_1=1}^k p_{m_1}\delta^{-\sigma-m_1+1} \;\mbox{ and }\;
t=2\sum_{m_2=1}^k q_{m_2}\delta^{-2\sigma-m_2+1},
$$
where $p_{m_1},q_{m_2}\in\mathbb{N}$. For such $s$, $t$ and $\xi$
we have
\begin{eqnarray*}
\begin{aligned}
s\xi-t\xi^2/2&=
\sum_{m_1=1}^k\sum_{r=1}^kp_{m_1}\ell_r\delta^{r-m_1}\\
&+ \sum_{m_1=1}^k\varepsilon p_{m_1}\delta^{-\sigma-m_1+1}-
\sum_{m_2=1}^k\sum_{r,r'=1}^k
q_{m_2}\ell_r\ell_{r'}\delta^{r+r'-m_2-1}\\
&-\sum_{m_2=1}^k\sum_{r=1}^k\varepsilon\ell_r
q_{m_2}\delta^{-\sigma+r-m_2} -\sum_{m_2=1}^k\varepsilon^2
q_{m_2}\delta^{-2\sigma-m_2+1}.
\end{aligned}
\end{eqnarray*}
In order for the above expression for the phase to belong to a
small (say $1/10$) neighbourhood of $\mathbb{N}$, it suffices for
each of the summands to either be integers, or be bounded in
absolute value by a sufficiently small constant (depending only on
$k$). As before, this places size restrictions on the integers
$p_{m_1}$ and $q_{m_2}$. It is here where we use the fact that
$1/\delta$ is an integer.

We now define
$$X_k=\{p_1\delta^{-\sigma}+\cdots+p_k\delta^{-\sigma-k+1}:p_1,\hdots,p_k\in\mathbb{N}\mbox{
with } p_1,\hdots,p_k\lesssim\delta^{\sigma-1}\},$$ and
$\Omega_k$ to be an $O(1)$-neighbourhood of
\begin{eqnarray*}
\begin{aligned}
\Lambda_k:=\{ (x,t)
: x\in X_k^{n-1},\;\;
t = 2(q_1\delta^{-2\sigma}+\cdots+q_k\delta^{-2\sigma-k+1});\;\;
 q_1,\hdots,q_k\lesssim\delta^{2\sigma -1}\}.
\end{aligned}
\end{eqnarray*}
Arguing as before we find that $|e^{it
\Delta}f_k(x)|\sim\delta^{k(n-1)(1-\sigma)}$ whenever
$(x,t)\in\Omega_k$.

We note a further ``self-similarity" in the family of sets
$\Lambda_k$. Observe that for $k\geq 2$, $\Lambda_k$ is a disjoint
union of rescaled and translated copies of $\Lambda_{k-1}$, where
$\Lambda_1$ coincides with the set $\Lambda$ defined previously.

The above examples may be generalised substantially by choosing
the family of functions $g_k$ to be self-similar with respect to
more general families of affine transformations that in
particular may depend on the index $k$.
\section{Weighted estimates for extension operators}
Let $n\geq 2$ and $S$ be a bounded hypersurface in $\mathbb{R}^n$
with everywhere non-vanishing Gaussian curvature (for instance,
$S$ could be the base of a paraboloid or a small portion of the unit sphere
$\mathbb{S}^{n-1}$). If we denote by $\sigma$ the induced Lebesgue
measure on $S$, then we may define the \emph{extension operator
associated to $S$} to be the mapping
$g\mapsto\widehat{gd\sigma}$ where
$$
\widehat{gd\sigma}(x)=\int g(\xi)e^{-2\pi
ix\cdot\xi}d\sigma(\xi),$$ $g\in L^{1}(S)$ and $x\in\mathbb{R}^n$.

\subsection{Rates of decay}
\vspace{0.5cm}
\noindent
There has recently been considerable interest in studying weighted
$L^{2}$ inequalities for the extension operator which take the general form
\begin{equation}\label{general}
\int_{\mathbb{B}}| \widehat{gd\sigma }(Rx)| ^{2}d\mu (x)\leq
\frac{C(\mu)}{R^{\gamma }}\left\| g\right\| _{L^{2}(S)}^{2}
\end{equation}
where $\mu$ is a positive measure supported on the unit ball
$\mathbb{B}$ of $\mathbb{R}^{n}$, $R\geq 1$, $\gamma$ is a
suitable rate-of-decay exponent and $C(\mu)$ is a constant
depending only on $\mu$. See for example \cite{Mattila},
\cite{Sj}, \cite{B}, \cite{BRV1}, \cite{CSo}, \cite{Wolff},
\cite{CSV}, \cite{BCSV}, \cite{E} and \cite{IR}. A specific
instance of this type of inequality, which is of particular
interest in geometric measure theory, concerns the relation
between the exponents $\gamma\geq 0$ and $0\leq \eta\leq n$ such
that for each $S$ there exists a constant $C$ (depending on $S$,
$\gamma$ and $\eta$)
for which
\begin{equation}\label{2.1}
\int_{\mathbb{B}} |\widehat{gd \sigma}(Rx) |^2 d \mu \leq
\frac{C}{ R^{\gamma}}\;\sup_{x\in\mathbb{R}^n, r>0 }
\left\{\frac{\mu (B(x,r))}{r^{ \eta }} \right\} \;||g||_2^2
\end{equation}
holds for all $g\in L^2(S)$, all $R\geq 1$ and all Borel measures $\mu$ supported
in $\mathbb{B}$. In particular, for each
$0\leq\eta\leq n$ it is of interest to determine the exponent
$\gamma(\eta)$ which is defined to be the supremum of the numbers $\gamma$
for which \eqref{2.1} holds for some constant $C$.

In two dimensions it is known that
\[
\gamma(\eta)=\left\{\begin{array}{lll}
\eta/2, &\;\;\;\; 1\leq\eta\leq 2\;\;\;\;\;\;\mbox{Wolff \cite{Wolff}}\\
1/2,& 1/2\leq\eta<1\;\;\;\;\;\; \mbox{Mattila \cite{M2}}  \\
\eta, &\;\;\;\; 0\leq\eta<1/2\;\; \mbox{Mattila \cite{M2}}.  \\
\end{array}
\right.
\]
We note that this piecewise linear function was originally
computed for $S=\mathbb{S}^{n-1}$. It is however implicit in the
arguments given in \cite{M2} and \cite{Wolff} that the same is true
for general $S$ of the type discussed here.

In higher dimensions it is well known that $\gamma(\eta)=\eta$ for
$0\leq\eta\leq\tfrac{n-1}{2}$ (Mattila \cite{M2}), and that
$\gamma(n)=n-1$ (Sj\"olin \cite{Sj}). However, for
$\tfrac{n-1}{2}<\eta<n$ the currently known upper and lower bounds
for $\gamma(\eta)$ do not coincide. Lower bounds in this region
were obtained by Sj\"olin \cite{Sj} and Bourgain \cite{B}, and those of
\cite{B} have been improved recently by Erdo$\tilde{\mbox{g}}$an \cite{E}.
Examples leading to upper bounds were also observed by Mattila \cite{M2},
Sj\"olin \cite{Sj}, Katz and Tao \cite{KT} and more
recently by Iosevich and Rudnev \cite{IR}. Our purpose here is to
improve these upper bounds further by using our special solutions
to the Schr\"{o}dinger equation from Section 2 by showing that the graph of
$\gamma(\eta)$ does not lie above the line segment joining the
points $(\tfrac{n-1}{2},\tfrac{n-1}{2})$ and $(n,n-1)$.

\begin{proposition}\label{upperbound} If for all bounded hypersurfaces $S$,
\eqref{2.1} holds for all $g\in L^2(S)$, all $R\geq 1$ and all Borel measures
$\mu$ supported in $\mathbb{B}$, and if $(n-1)/2 < \eta < n$, then
$$
\gamma \leq (\eta+1)\left(\frac{n-1}{n+1}\right).
$$
\end{proposition}

It is important to point out that we give upper bounds for the
general problem \eqref{2.1} by construction of examples on a
particular hypersurface, the paraboloid. These examples do not
appear to extend in a routine manner to general hypersurfaces, or
even to the specific case of the sphere $S=\mathbb{S}^{n-1}$, due
to number-theoretic issues. The case of the sphere is the
principal interest of \cite{IR} as it relates to the classical
distance set conjecture of Falconer \cite{Fal}.

\begin{proof}
We take $S$ to be the section of the
paraboloid
\begin{equation}\label{paraboloid}
\{\xi=(\xi',\xi_n)\in\mathbb{R}^{n-1}\times\mathbb{R}:\xi_n=|\xi'|^2/2,
\; 0\leq \xi_1,\hdots\xi_{n-1}\leq 1\},\end{equation} and observe
that
\begin{equation}\label{red}
\widehat{gd\sigma}(x)= \int_{|\xi'|\leq 1}e^{-2\pi
ix'\cdot\xi'+\pi i|\xi'|^{2}x_{n}}
\hat{f}(\xi')d\xi'=e^{ix_n\Delta}f(x'),
\end{equation}
where $x=(x',x_n)\in\mathbb{R}^{n-1}\times\mathbb{R}$ and
$\hat{f}(\xi')=g(\xi',|\xi'|^{2})(1+|\xi'|^{2})^{1/2}$. Now by
Plancherel's Theorem, inequality \eqref{2.1} may be written as
\begin{equation}\label{pa}
\int_{\mathbb{B}} |e^{iRx_n\Delta}f(Rx')|^2 d \mu \leq \frac{C}{
R^{\gamma}}\;\sup_{x\in\mathbb{R}^n, r>0 } \left\{\frac{\mu
(B(x,r))}{r^{ \eta }} \right\} \;||f||_2^2,
\end{equation}
and so we may test \eqref{pa} on functions $f$ of the form
\eqref{exs} with $\delta=1/R$. For such a function
$|e^{iRx_n\Delta}f(Rx')|\gtrsim R^{-(n-1)(1-\sigma)}$ for all
$x\in\widetilde{\Omega}$, where $\widetilde{\Omega}$ is an
$O(1/R)$ neighbourood of
$$
\{x\in \mathbb{R}^{n}: x'\in \widetilde{X}^{n-1}\mbox{ and } x_n =
2q R^{2 \sigma-1} \mbox{ where }q \in \mathbb{N} \mbox{ and }
q\lesssim R^{1-2\sigma}\}
$$
and
$$
\widetilde{X}=\{pR^{\sigma-1}:p\in\mathbb{N}\mbox{ with }
p\lesssim R^{1-\sigma}\}.$$ Notice that
$R\widetilde{\Omega}=\Omega$.

Setting $d\mu(x)=\chi_{\widetilde{\Omega}}(x)dx$ we therefore
obtain
$$\int_{\mathbb{B}} |e^{iRx_n\Delta}f(Rx')|^2 d \mu\gtrsim
R^{-2(n-1)(1-\sigma)}|\widetilde{\Omega}|\sim
R^{-2(n-1)(1-\sigma)-\sigma(n+1)}.
$$
Furthermore,
\[
\sup_{x\in\mathbb{R}^n, r>0 } \left\{\frac{\mu (B(x,r))}{r^{ \eta
}} \right\}\sim\left\{\begin{array}{ll}
R^{-\sigma(n+1)}, & \frac{n-1}{2}\leq\eta\leq n-\sigma(n+1)\\
R^{\eta-n},& n-\sigma(n+1)\leq\eta\leq n \\
\end{array}
\right.
\]
as the reader will easily verify.

Using these calculations and the fact that $\|f\|_2^2\sim
R^{-(n-1)(1-\sigma)}$, we may deduce that a necessary condition
for \eqref{pa} (and thus \eqref{2.1}) to hold for all $R\geq 1$ is
that
$$
\gamma\leq (\eta+1)\left(\frac{n-1}{n+1}\right),
$$
as claimed.
\end{proof}

It is conceivable that the upper bound of Proposition \ref{upperbound}
may be improved by considering the more sophisticated special solutions of
the Schr\"{o}dinger equation \eqref{ecuacion} described at the end
of Section 2, but we do not pursue this point further.

We note that the examples of Section 2 furnish explicit necessary conditions
on the functional $C(\mu)$ so that an inequality of the form \eqref{general}
might hold for some given $R$. These conditions are not simply upper bounds
on how much mass $\mu$ can put on various eccentric tubes, but instead on
how much mass can be put on various {\em{arrangements}} of eccentric tubes
dictated by the set $\Omega$ and its variants. On the other hand, if we
take for example $\gamma = n-1$, and we demand validity of \eqref{general}
for {\em{all}} $R \geq 1$ then the examples we present here offer no further
necessary conditions than do the ``simple" examples where the testing
function is essentially the characteristic function of a {\em{single}}
product of intervals.
\subsection{Morrey--Campanato weights}
\label{restcompact}

The Stein--Tomas restriction  theorem (see \cite{St}), in its equivalent dual form,
states that for $r\geq 2(n+1)/(n-1)$, there is a constant $C$ for which
\begin{equation}\label{ST}
\|\widehat{g d\sigma}\|_{L^r(\mathbb R^n)}
\leq C\ \|g\|_{L^2 (S)},
\end{equation}
for all $g\in L^2 (S)$. This  estimate  may  be viewed using duality as  the
weighted  estimate
\begin{equation}\label{WST}
\int_{\mathbb{R}^n} | \widehat{g d\sigma}(x)|^2 V(x)dx
\leq C\ \|V\|_{L^{p} ({\mathbb R^n})} \int_S |g |^2 d \sigma,
\end{equation}
whenever $1 \leq p \leq  (n+1)/2$ and $V\in  L^p$.

In \cite{RV91} Ruiz  and Vega considered (see also \cite{cha+s},
\cite{chi+r} and \cite{RV94} for $\alpha=2$) extending \eqref{WST}
by replacing the $L^p$  norm  of $V$ by certain
\emph{Morrey--Campanato}  norms. These norms play a role in the
theory of unique continuation -- see \cite{Ke} and \cite{Wolff2}.
The Morrey--Campanato classes, which are denoted by $\mathcal{L}^{\alpha,p}$,
for $\alpha>0$ and $1\le p\le n/\alpha$,
are given  by
$$
\mathcal{L}^{\alpha,p} = \left\{ V\in L_{\ell oc}^p(\mathbb{R}^n)
:\|V\|_{\mathcal{L}^{\alpha,p}}<\infty\right\},
$$
where
$$
\|V\|_{\mathcal{L}^{\alpha,p}} = \sup_{x\in\mathbb{R}^n,r>0}\;
r^{\alpha}\left(r^{-n}\int_{B(x,r)}|V(y)|^p\ dy\right)^{1/p}.
$$
Notice that
$\mathcal{L}^{\alpha,n/\alpha}=L^{n/\alpha}(\mathbb{R}^n).$ We
also remark  that  for  $p< n/\alpha$  the  class
$\mathcal{L}^{\alpha,p}$ contains  the  Lorentz space
$L^{n/\alpha, \infty}(\mathbb{R}^n)$. Furthermore, when $p=1$, it is natural
to think of $V$ as a measure $\mu$ rather than a locally $L^1$ function,
and then $$\|\mu\|_{\mathcal{L}^{\alpha,1}} = \sup_{x\in\mathbb{R}^n,r>0}\;
\left\{\frac{\mu(B(x,r)}{r^{n-\alpha}}\right\},$$ (similar to what was
considered in the previous subsection).

\begin{proposition}\label{RV}
If $n \geq 2$, $\tfrac{2n}{n+1} < \alpha \leq \displaystyle n$, $p
\geq 1$ and $\displaystyle \tfrac{\alpha}{n}\leq \tfrac{1}{p} <
\tfrac{2(\alpha-1)}{n-1}$, there exists a constant $C$ (depending
on $S$, $p$ and $\alpha$, but independent of $g$ and $V$) such
that
\begin{equation}\label{extension}
\int_{\mathbb{R}^n} |\widehat{g d\sigma}(x)|^2 V(x) dx \leq C\
\|V\|_{\mathcal{L}^{\alpha,p}}\int_S |g |^2.
\end{equation}
\end{proposition}



The proof given by Ruiz and Vega in \cite{RV91} was for the case
$S = \mathbb{S}^{n-1}$, but all the estimates used go through in
the more general case of nonvanishing Gaussian curvature. There are
three ``endpoint" cases $(\alpha, \tfrac{1}{p}) = (\tfrac{2n}{n+1},
\tfrac{2}{n+1}), \; (\tfrac{n+1}{2}, 1)$ and $ (n,1)$ respectively.
The first of these is the Stein--Tomas restriction theorem \eqref{WST},
the second corresponds to a rescaled version of Mattila's result
(\cite{M2}) that $\gamma(\eta)=\eta$ for $\eta = \tfrac{n-1}{2}$ which
was discussed in the previous subsection, and the third is trivial. (The key
difficulty to be overcome in \cite{RV91} was the failure of the
Morrey--Campanato spaces to interpolate nicely; see also \cite{BRV}.)

We now consider the sharpness of the condition $\displaystyle
\tfrac{1}{p} < \tfrac{2(\alpha-1)}{n-1}$ in Proposition \ref{RV}.
(See \cite {Ke} and \cite{Wolff2} for related issues concerning Carleman
estimates.)

In the first place, a straightforward modification of the standard
``Knapp" counterexample provides the following necessary condition
for $\alpha < 2.$ This condition gives optimal results in
dimensions $n=2,3$ with the possible exception of the line
$\tfrac{1}{p}=\tfrac{2(\alpha-1)}{n-1}$.
\begin{lemma}
\label{nrestriccion} Let $n\geq 2$, $\alpha < 2 $ and suppose that
\eqref{extension} holds for some $S$. Then $\displaystyle
\tfrac{1}{p}\leq \tfrac{2(\alpha-1)}{n-1}.$
\end{lemma}

\begin{proof}
Let $g$ be the characteristic function of a $\delta$-cell on $S$;
then $\|g\|_2^2 \sim \delta^{n-1}$, while $|\widehat{g d \sigma}(x)|
\geq C \delta^{n-1}$ on a tube of sides $\delta^{-1} \times \delta^{-1}
\cdots \times \delta^{-1} \times \delta^{-2}$. Take $V$ to be the
characteristic function of this tube. Then
$$
\displaystyle
\|V\|_{\mathcal{L}^{\alpha,p}}
\sim
\left\{
\begin{array}{ll}
\displaystyle
\delta^{-\alpha}&{\rm if}\ \frac{1}{p}\ge
\frac{\alpha}{n-1},
\\[2ex]
\displaystyle
\delta^{-2\alpha+\frac{n-1}{p}}&{\rm if}\ \frac{1}{p}\le \frac{\alpha}{n-1}.
\end{array}
\right.
$$


The claim follows from the second of these estimates upon taking
$\delta$ small.
\end{proof}


These standard examples only have significance when $\alpha  < 2$.
As we shall now see, for  $\alpha \geq 2$ and $n\geq 4$, the
special solutions introduced in Section \ref{especiales} provide
us with further necessary conditions.  As in the previous
subsection, this will be achieved by taking $S$ to be a bounded
subset of the paraboloid. (Again, we point out that these examples
do not appear to extend in a routine manner to other curved
submanifolds of the type we consider.)

\begin{proposition}
\label{paraboloide} Suppose that $n\geq 4$  and  $S$  is the
section of the paraboloid given by \eqref{paraboloid}. If $ \alpha
\geq 2$  and \eqref{extension} holds, then $\displaystyle
\tfrac{1}{p} \leq \tfrac{2\alpha}{n+1}.$
\end{proposition}
\begin{proof}
As we did in \eqref{red}, we can write
$\widehat{gd\sigma}(x)=
e^{ix_n\Delta}f(x'),$ where
$x=(x',x_n)\in\mathbb{R}^{n-1}\times\mathbb{R}$ and
$\widehat{f}(\xi')=g(\xi',|\xi'|^{2})(1+|\xi'|^{2})^{1/2}$ and
therefore, inequality \eqref{extension} may be written as
\begin{equation}
\label{trozo}
\|e^{ix_n\Delta}f(x')\|_{L^2(V)} \le C\
\|V\|^{1/2}_{\mathcal{L}^{\alpha,p}} \|\widehat{f}\|_{L^2(\mathbb{R}^{n-1})}.
\end{equation}
Taking $f$ as in \eqref{exs} and $V$ as the characteristic
function of the set $\Omega$ defined in Section \eqref{especiales}
(see Figure \ref{Omega}), we have that
\begin{equation}
\label{trozoleft}
\|e^{ix_n\Delta}f(x')\|_{L^2(V)}
\sim \delta^{(n-1)(1-\sigma)}|\Omega|^{1/2} =
\delta^{\frac{(1-\sigma)(n-1)}{2}+\sigma-\frac{1}{2}}.
\end{equation}
On the other hand,
\begin{equation}
\label{trozoright}
\|\widehat{f}\|_{L^2(\mathbb{R}^{n-1})}
\sim
\delta^{\frac{(1-\sigma)(n-1)}{2}}.
\end{equation}
For $p\le n/\alpha$, $0<\sigma<1/2,$ and $\delta$ small
and positive we have
\begin{equation}
\label{morreyxt}
\|V\|_{\mathcal{L}^{\alpha,p}}
\sim
\max\{1,\ \delta^{\frac{\sigma(n+1)}{p}-\alpha}\}.
\end{equation}
From \eqref{trozoleft}, \eqref{trozoright} and
\eqref{morreyxt}, and the fact that $0<\sigma<1/2,$ we see that
a necessary condition for \eqref{trozo} to hold is
$\displaystyle \tfrac{1}{p} \leq \tfrac{2\alpha}{n+1}.$
\end{proof}

\begin{figure}[h]
{\centerline{\psfig{figure=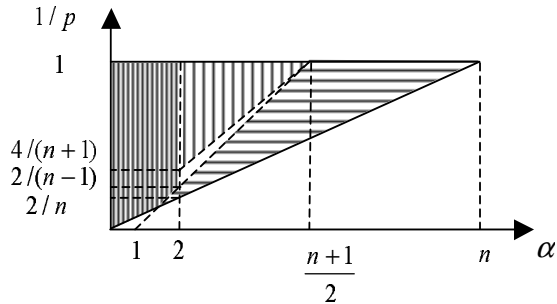,height=4cm,width=7.3cm}}}
\caption{Let $S$ be the section of the paraboloid defined in (\ref{paraboloid}).
The regions with horizontal lines and vertical lines correspond to the cases
where estimate \eqref{extension} is true and false respectively.}
\label{ptrozo}
\end{figure}

Figure \ref{ptrozo} shows the positive and the negative results
discussed here when $ \alpha \leq \tfrac{n}{p}$ and $S$ is the
section of the paraboloid \eqref{paraboloid}.
We shall address the analogous questions for the full paraboloid -- i.e.
the initial value problem for the Schr\"odinger equation \eqref{ecuacion}
where the initial data $f$ is not assumed to have Fourier transform with
compact support -- in a forthcoming paper. As a consequence of those results
we shall obtain some further necessary conditions for \eqref{extension} to
hold in the case of the sphere $\mathbb{S}^{n-1}$.
\bibliographystyle{plain}

\end{document}